\documentclass[10pt,twocolumn,twoside]{IEEEtran}

\usepackage{hyperref} 
\usepackage{amsmath,amssymb,amsfonts}%
\usepackage{theorem}%
\usepackage[dvipsnames]{xcolor}
\usepackage{graphicx}%

\newtheorem{definition}{Definition}[section]%
\newtheorem{theorem}[definition]{Theorem}%
\newtheorem{proposition}[definition]{Proposition}%
\newtheorem{lemma}[definition]{Lemma}%
\newtheorem{assumption}[definition]{Assumption}%
\newtheorem{corollary}[definition]{Corollary}%
{\theorembodyfont{\rmfamily} \newtheorem{remark}[definition]{Remark}}%
{\theorembodyfont{\rmfamily} }%

\newcommand{\trn}{^{\scriptscriptstyle \top}}%
\newcommand{\esssup}{\operatorname*{ess\;sup}}%
\newcommand{\inner}{\mathrm{int}}%
\newcommand{\tm}{\times}%
\newcommand{\N}{\mathbb{N}}%
\newcommand{\R}{\mathbb{R}}%

\newcommand{\KC}{\mathcal{K}}%
\newcommand{\EC}{\mathcal{E}}%
\newcommand{\VC}{\mathcal{V}}%
\newcommand{\NC}{\mathcal{N}}%
\newcommand{\LC}{\mathcal{L}}%
\newcommand{\UC}{\mathcal{U}}%

\newcommand{\rmd}{\mathrm{d}}%
\newcommand{\rmD}{\mathrm{D}}%
\newcommand{\rme}{\mathrm{e}}%
\newcommand{\diag}{\mathrm{diag}}%
\newcommand{\A}{\mathcal{A}}%
\newcommand{\abs}[1]{\left|#1\right|}%

\let\ol=\overline%
\let\ul=\underline%
\newcommand{\Rn}[1][n]{\R^{#1}}%
\def\K{\mathcal{K}}%
\def\KL{\mathcal{KL}}%
\def\Kinf{\mathcal{K}_\infty}%

\allowdisplaybreaks%
\begin{document}

\title{Small-gain theorem for stability, cooperative control and distributed observation of infinite networks}

\author{Navid~Noroozi, Andrii~Mironchenko, Christoph~Kawan, Majid~Zamani,~\IEEEmembership{Senior Member,~IEEE}

\thanks{N.~Noroozi is with the Otto-von-Guericke University Magdeburg, Laboratory for Systems Theory and Automatic Control, Magdeburg, Germany; email: \texttt{navid.noroozi@ovgu.de}. His work is supported by the DFG through the grant WI 1458/16-1.}
\thanks{A.~Mironchenko is with Faculty of Computer Science and Mathematics, University of Passau, 94032 Passau, Germany; email: \texttt{andrii.mironchenko@uni-passau.de}. His work is supported by the DFG through the grant MI 1886/2-1.}
\thanks{C.~Kawan is with the Institute of Informatics, Ludwig Maximilian University of Munich, Germany; email: \texttt{christoph.kawan@lmu.de}. His work is supported by the DFG through the grant ZA 873/4-1.}%
\thanks{M.~Zamani is with the Computer Science Department, University of Colorado Boulder, CO 80309, USA. 
M.~Zamani is also with the Computer Science Department, Ludwig Maximilian University of Munich, Germany; email: {\tt\small majid.zamani@colorado.edu}. His work is supported in part by the DFG through the grant ZA 873/4-1 and the H2020 ERC Starting Grant AutoCPS (grant agreement No. 804639).}}

\date{}%
\maketitle%

\begin{abstract}
Motivated by a paradigm shift towards a \emph{hyper-connected} world, we develop a \emph{computationally tractable} small-gain theorem for a network of \emph{infinitely} many systems, termed as infinite networks. The proposed small-gain theorem addresses exponential input-to-state stability with respect to closed sets, which enables us to analyze diverse stability problems in a unified manner. The small-gain condition, expressed in terms of the spectral radius of a gain operator collecting all the information about the internal Lyapunov gains, can be numerically computed for a large class of systems in an efficient way. To demonstrate broad applicability of our small-gain theorem, we apply it to the stability analysis of infinite \emph{time-varying} networks, to consensus in \emph{infinite}-agent systems, as well as to the design of \emph{distributed observers} for infinite networks.
\end{abstract}

\section{Introduction}%

Emerging technologies such as the Internet of Things, Cloud computing, 5G communication, and so on are expected to encompass almost every aspect of our lives. Those advances will result in a paradigm shift towards a \emph{hyper-connected} world composed of a large number of smart networked systems providing us with much more \emph{autonomy} and \emph{flexibility}.
However, these benefits are obtained at the price of increasing \emph{complexity} and \emph{uncertainty}.
Examples of such smart networked systems include smart grids, connected vehicles, swarm robotics, and smart cities in which the participating agents may be plugged into and out from the network at any time. Thus, the size of such very large networks is unknown and possibly time-varying.

Most of these smart applications are \emph{safety-critical}. This calls for a rigorous analysis and synthesis of such networks of systems. However, standard tools for stability analysis/stabilization of control systems do not scale well to these large-scale complex systems~\cite{Sarkar.2018,Bamieh.2012,Jovanovic.2005b}.
A promising approach to address this critical issue is to over-approximate a finite but very large network by an \emph{infinite network}, and then control this over-approximated system; see e.g.~\cite{CIZ09,Jovanovic.2005b,DAndrea.2003}.

Current results on the stability analysis and control of infinite networks are mostly concerned with spatially invariant and/or linear systems~\cite{CIZ09,DAndrea.2003}.
A striking progress in the infinite-dimensional input-to-state stability (ISS) theory~\cite{DaM13,MiW18b,KaK16b,JNP18} (see~\cite{MiP19} for a recent survey on this topic) blended with the powerful nonlinear small-gain criteria for stability analysis of finite networks of nonlinear systems \cite{JMW96,DRW07,Dashkovskiy.2010} create a foundation for the development of stability conditions for infinite networks of general nature without assuming linearity and/or spatial invariance of the systems.

The case of infinite networks is much more complex, as the gain operator, collecting the information about the internal gains, acts in an infinite-dimensional space, in contrast to couplings of just $N\in\N$ systems of arbitrary nature (possibly infinite-dimensional).
This calls for a careful choice of the infinite-dimensional state space of the overall network, and motivates the use of the theory of positive operators on ordered Banach spaces for the small-gain analysis.

In~\cite{DaP19} it is shown that a countably infinite network of continuous-time input-to-state stable systems is ISS, provided that the gain functions capturing the influence of subsystems at each other are all less than identity, which is a very conservative condition.
In~\cite{DMS19a} it was shown that classic max-form strong small-gain conditions (SGCs) developed for finite networks in \cite{Dashkovskiy.2010} do not ensure stability of infinite networks, even for linear ones. To address this issue, more restrictive robust strong SGCs are developed in~\cite{DMS19a}. The small-gain theorems in~\cite{DaP19,DMS19a} are formulated in terms of ISS Lyapunov functions and a trajectory-based small-gain theorem for infinite networks is provided in~\cite{Mir19d}.

By contrast, for networks consisting of exponentially ISS systems, possessing exponential ISS Lyapunov functions with linear gains, it was shown in \cite{KMS19} that if the spectral radius of the gain operator is less than one, then the whole network is exponentially ISS and there is a coercive exponential ISS Lyapunov function for the whole network. This result provides a complete and \emph{nontrivial} generalization of \cite[Prop.~3.3]{Dashkovskiy.2011b} from finite networks to infinite ones. It deeply relies on the spectral theory of positive operators~\cite{K59}. The effectiveness of the main result in~\cite{KMS19} has been demonstrated by applications to nonlinear spatially invariant systems with sector nonlinearities and to the stability analysis of a road traffic network.

All of the above small-gain theorems for infinite networks address ISS with respect to the origin. A more general notion of input-to-state stability with respect to a closed set covers several further stability notions such as incremental stability, robust consensus/synchronization, ISS of time-varying systems as well as variants of input-to-output stability in a unified and generalized manner~\cite{Noroozi.2018a}. In this paper, we extend the main result of our recent work~\cite{KMS19} to ISS of infinite networks with respect to closed sets. This generalization widely extends the applicability of the small-gain result to several control theoretic problems including the stability analysis of infinite time-varying networks, consensus of infinite multi-agent systems, as well as the design of distributed observers for infinite networks which all are demonstrated in this work.

In the literature, the stability theory for nonlinear time-invariant systems and nonlinear time-varying systems is often presented separately.
Moreover, existing results on infinite networks are developed for time-invariant systems, although, practically speaking, time-variance is a more realistic assumption.
In this paper, as the first application of our result, we address exponential ISS for time-invariant and time-varying infinite networks within a \emph{unified} framework.

Distributed cooperative control has broad applications in various areas such as wireless sensor networks, mobile robots, power networks, social networks, etc~\cite{Kia.2019,Ren.2007}.
The current literature on distributed cooperative control mainly focuses on the case of networks of fixed number of agents, i.e. the size of physical networks can increase or decrease over time.
In several applications such as social networks, the size of the network, however, is time-varying, often huge and uncertain.
Such networks are known as ``\emph{open multi-agent systems}''~\cite{Varma.2018}.
To address the scalability issue of open multi-agent systems, here we overapproximate the network with a \emph{time-varying} number of components by countably infinite networks and then seek for consensus of the infinite number of agents.
In particular, we formulate a weighted average consensus problem~\cite{Ren.2007} for the infinite network as a stabilization problem with respect to a closed set, which is tractable by means of the small-gain approach developed in this work. The study of the error between the real open multi-agent system and its infinite approximation is however beyond the scope of this paper and is an interesting open problem itself.

Motivated by applications in surveillance and monitoring for spatially distributed systems such as environmental and agricultural monitoring, healthcare monitoring, and pollution source localization, for the third application of our result, we provide a methodology to address scalability issues in distributed estimation problems.
We assume that each subsystem has a local observer asymptotically converging to the true state of each subsystem, given perfect knowledge of the true states of neighboring subsystems. Formulating the state estimation as a stabilization problem with respect to a certain closed set, we show that if the couplings between subsystems are small enough, which is quantitatively expressed by our small-gain condition, then the state estimation problem can be solved.


This paper is organized as follows: First, relevant notation, discussions on well-posedness of infinite networks and an appropriate distance function with respect to a closed set in an infinite-dimensional state space are given in Section~\ref{sec_prelim}. The notion of exponential ISS with respect to a closed set for infinite-dimensional systems and related Lyapunov properties are presented in Section~\ref{sec:eISS}.
Technical results on the gain operator are made precise in Section~\ref{sec:The-gain-operator-and-its-properties}. In Section~\ref{sec:Small-gain-theorem}, the main result of the paper is presented. The effectiveness of our result is verified through applications to time-varying infinite networks, consensus problems and distributed observers in Section~\ref{sec:Applications}. In Section~\ref{sec:Conclusions}, we conclude the paper. 
The proofs of auxiliary results are provided in Appendix.

\section{Preliminaries}\label{sec_prelim}

\subsection{Notation}%

We write $\N = \{1,2,3,\ldots\}$ for the set of positive integers, $\R$ denotes the reals and $\R_+ := \{t \in \R : t \geq 0\}$ the nonnegative reals. For vector norms on finite- and infinite-dimensional vector spaces, we write $|\cdot|$. For associated operator norms, we use the notation $\|\cdot\|$. We write $A\trn$ for the transpose of a matrix $A$ (which can be finite or infinite). We typically use Greek letters for infinite matrices and Latin ones for finite matrices. Elements of $\R^n$ are by default regarded as column vectors and we write $x\trn \cdot y$ for the Euclidean inner product of two vectors $x,y \in \R^n$. We use the same notation for dot products of vectors with infinitely many components. By $\ell^p$, $p\in[1,\infty]$, we denote the Banach space of all real sequences $x = (x_i)_{i\in\N}$ with finite $\ell^p$-norm $|x|_p<\infty$, where $|x|_p = (\sum_{i=1}^{\infty}|x_i|^p)^{1/p}$ for $p < \infty$ and $|x|_{\infty} = \sup_{i\in\N}|x_i|$. We write $\ell^p_+ := \{ x = (x_i)_{i\in\N} \in \ell^p : x_i \geq 0,\ \forall i \in \N \}$.%

A more general class of $\ell^p$-spaces is defined as follows. Let $p\in[1,\infty)$, let $(n_i)_{i\in\N}$ be a sequence of positive integers and fix a norm $|\cdot|_i$ on $\R^{n_i}$ for every $i\in\N$. Then%
\begin{equation*}
  \ell^p(\N,(n_i)) := \Bigl\{ x = (x_i)_{i\in\N} : x_i \in \R^{n_i},\ \sum_{i=1}^{\infty}|x_i|_i^p < \infty \Bigr\}%
\end{equation*}
equipped with the norm%
$
  |x|_p := \Bigl(\sum_{i=1}^{\infty} |x_i|_i^p\Bigr)^{\frac{1}{p}}%
$
is a separable Banach space (can be proved using standard arguments see e.g.~\cite{DunSch57}).
Usually, we drop the index $i$ from the norm. If all $n_i$ are identical, say $n_i \equiv n$, we also write $\ell^p(\N,n)$.
Similarly, $\ell^{\infty}(\N,(n_i))$ can be defined.%

We write $L^{\infty}(\R_+,\R^n)$ for the Banach space of essentially bounded measurable functions from $\R_+$ to $\R^n$. If $X$ is a Banach space, we write $r(T)$ for the spectral radius of a bounded linear operator $T:X \rightarrow X$ and $L(X)$ for the space of all bounded linear operators on $X$.
The notation $C^0(X,Y)$ stands for the set of all continuous mappings $f:X \rightarrow Y$ between metric spaces $X$ and $Y$.
Given a metric space $X$, we write $\inner\, A$ for the interior of a subset $A \subset X$.
The right upper (resp.\ lower) Dini derivative of a function $\gamma:\R \rightarrow \R$ at $t\in\R$ is denoted by $\rmD^+\gamma(t)$ (resp.\ $\rmD_+\gamma(t)$); see~\cite{KMS19} for their definitions.
We will consider $\K, \Kinf$, and $\KL$ comparison functions, see~\cite[Chapter 4.4]{Khalil.2002} for definitions.

\subsection{Infinite interconnections}%

We study interconnections of countably many systems, each given by a finite-dimensional ordinary differential equation (ODE). Using $\N$ as the index set (by default), the $i$th subsystem is written as%
\begin{equation}\label{eq_ith_subsystem}
  \Sigma_i:\quad \dot{x}_i = f_i(x_i,\bar{x},u_i).%
\end{equation}
The family $(\Sigma_i)_{i\in\N}$ comes together with a number $p \in [1,\infty]$ and sequences $(n_i)_{i\in\N}$, $(m_i)_{i\in\N}$ of positive integers so that the following assumptions hold with $X := \ell^p(\N,(n_i))$ for a specified sequence of norms on the spaces $\R^{n_i}$:%
\begin{itemize}
\item The state vector $x_i$ of $\Sigma_i$ is an element of $\R^{n_i}$.%
\item The internal input vector $\bar{x}$ is an element of $X$.%
\item The external input vector $u_i$ is an element of $\R^{m_i}$.%
\item The right-hand side $f_i:\R^{n_i} \tm X \tm \R^{m_i} \rightarrow \R^{n_i}$ is a continuous function.%
\item Unique local solutions of the ODE \eqref{eq_ith_subsystem} exist for all initial states $x_{i0} \in \R^{n_i}$ and all continuous $\bar{x}(\cdot)$ and locally essentially bounded $u_i(\cdot)$ (which are regarded as time-dependent inputs). We denote the corresponding solution by $\phi_i(\cdot,x_{i0},(\bar{x},u_i))$.%
\end{itemize}

The values of the function $f_i$ can be independent of certain components of the input vector $\bar{x}$. We write $I_i$ for the set of indices $j\in\N$ so that $f_i(x_i,\bar{x},u_i)$ is non-constant with respect to the component $x_j$ of $\bar{x}$, and without loss of generality we assume that $i \notin I_i$ (note that $f_i$ depends on $x_i$ explicitly).%

In the ODE~\eqref{eq_ith_subsystem}, we consider $\bar{x}(\cdot)$ as an \emph{internal input} and $u_i(\cdot)$ as an \emph{external input} (which may be a disturbance or a control input). The interpretation is that the subsystem $\Sigma_i$ is affected by a certain set of neighbors, indexed by $I_i$, and its external input. We note that the set $I_i$ does \emph{not} have to be finite, implying that subsystem $i$ can be connected to infinitely many other subsystems.%

To define the interconnection of the subsystems $\Sigma_i$, we consider the state vector $x = (x_i)_{i\in\N} \in X = \ell^p(\N,(n_i))$, the input vector $u = (u_i)_{i\in\N} \in \ell^q(\N,(m_i))$ for some $q \in [1,\infty]$ and the right-hand side $f(x,u) := (f_1(x_1,\bar{x},u_1),f_2(x_2,\bar{x},u_2),\ldots)$. The interconnection is then written as%
\begin{equation}\label{eq_interconnection}
  \Sigma:\quad \dot{x} = f(x,u).%
\end{equation}
The class of admissible control functions is defined as%
\begin{align}\label{eq:Input-space}
  \UC := \bigl\{ u:\R_+ \rightarrow U : &u \mbox{ is strongly measurable} \nonumber\\
  &\qquad\mbox{ and essentially bounded} \bigr\},%
\end{align}
and we equip this space with the sup-norm%
\begin{equation*}
  |u|_{q,\infty} := \esssup_{t\geq0}|u(t)|_q.%
\end{equation*}
  
A continuous mapping $\xi:I \rightarrow X$, defined on an interval $I = [0,T_*)$ with $T_* \in (0,\infty]$, is called a \emph{solution} of the infinite-dimensional ODE \eqref{eq_interconnection} with initial value $x^0 \in X$ for the external input $u\in\UC$ provided that the two conditions%
\begin{equation*}
  f(\xi(t),u(t)) \in X \mbox{\ \ and\ \ } \xi(t) = x^0 + \int_0^t f(\xi(s),u(s)) \rmd s%
\end{equation*}
hold for all $t \in I$, where the integral is the Bochner integral for Banach space valued functions. For the theory of the Bochner integral, a reader may consult, e.g., \cite{ABH11}.%

If for each $x^0 \in X$ and $u \in \UC$ a unique (local) solution exists, we say that the system is \emph{well-posed} and write $\phi(\cdot,x^0,u)$ for any such solution. As usual, we consider the maximal extension of $\phi(\cdot,x^0,u)$ and write $I_{\max}(x^0,u)$ for its interval of existence. We say that the system is \emph{forward complete} if $I_{\max}(x^0,u) = \R_+$ for all $(x^0,u) \in X \tm \UC$.%

We note that~\cite[Thm.~3.2]{KMS19} provides sufficient conditions for well-posedness of $\Sigma$.%

%

\subsection{Distances in sequence spaces}%

Let $X=\ell^p(\N,(n_i))$ for a certain $p\in[1,\infty)$. Consider nonempty closed sets $\A_i \subset \R^{n_i}$, $i\in\N$. For each $x_i \in \R^{n_i}$ we define the distance of $x_i$ to the set $\A_i$ by%
\begin{equation*}
  |x_i|_{\A_i} := \inf_{y_i\in \A_i}|x_i - y_i|.
\end{equation*}
Now we define the set%
\begin{eqnarray}\label{eq:Overall_set}
  \A := \{x \in X: \ x_i \in \A_i,\ i\in\N\} =  X \cap (\A_1\times \A_2 \times \ldots).%
\end{eqnarray}
If $\A \neq \emptyset$, we can define the distance from any $x\in X$ to $\A$ as%
\begin{eqnarray}\label{eq:Distance_to_overall_set}
  |x|_{\A} := \inf_{y\in \A}|x-y|_p = \inf_{y\in \A} \Big(\sum_{i=1}^\infty |x_i-y_i|^p \Big)^{\frac{1}{p}}.%
\end{eqnarray}

\begin{lemma}\label{lem:Alternative-A-representation} 
Let $X=\ell^p(\N,(n_i))$ for a certain $p\in[1,\infty)$.
Assume that $\A$ defined by \eqref{eq:Overall_set} is nonempty. Then for any $x \in X$
\begin{eqnarray}\label{eq:Distance_to_overall_set_Formula}
  |x|_{\A} = \Big(\sum_{i=1}^\infty |x_i|_{\A_i}^p \Big)^{\frac{1}{p}} < \infty.
\end{eqnarray}
\end{lemma}


Note that if $\A=\{0\}$, then we have $|x|_{\{0\}} = |x|_p$.%

\section{Exponential input-to-state stability}\label{sec:eISS}

Having a well-posed interconnection \eqref{eq_interconnection} with state space $X = \ell^p(\N,(n_i))$ and external input space $U = \ell^q(\N,(m_i))$ for $p,q \in [1,\infty)$, we aim to study the stability of the interconnected system with respect to a closed set $\A \subset X$. For this purpose, we introduce the notions of input-to-state stability and exponential input-to-state stability with respect to a set $\A$.%


\begin{definition}
Given a nonempty closed set $\A \subset X$, the system $\Sigma$ is called%
\begin{itemize}
\item \emph{input-to-state stable (ISS) w.r.t.~$\A$} if it is forward complete and there are functions $\beta \in \KC\LC$ and $\gamma \in \KC$ such that for any initial state $x^0 \in X$ and any $u\in\UC$ the corresponding solution satisfies%
\begin{equation*}
  |\phi(t,x^0,u)|_\A \leq \beta(|x^0|_\A,t) + \gamma(|u|_{q,\infty}) \mbox{\quad for all\ } t \geq 0.%
\end{equation*}
\item \emph{exponentially input-to-state stable (eISS) w.r.t.~$\A$} if it is ISS w.r.t.~$\A$ with a $\KC\LC$-function $\beta$ of the form $\beta(t,r) = M\rme^{-at}r$ for some $a,M>0$.%
\end{itemize}
\end{definition}

\begin{remark}
\label{rem:Exponential ISS and practical ISS} 
Let $\A$ be a bounded and closed set. Define $\|\A\|:=\sup_{x\in\A} |x|_p$ and note that for all $x \in X$ the inequality%
\begin{eqnarray}
|x|_p - \|\A\| \leq |x|_{\A} \leq |x|_p + \|\A\|
\label{eq:Relations_between_norms_1}
\end{eqnarray}    
holds. Suppose that \eqref{eq_interconnection} is eISS w.r.t. $\A$. Then the estimate%
\begin{equation*}
  |\phi(t,x^0,u)|_\A \leq M\rme^{-at}|x^0|_\A + \gamma(|u|_{q,\infty}), \quad t \geq 0%
\end{equation*}
implies that for $t \geq 0$
\begin{equation*}
  |\phi(t,x^0,u)|_p - \|\A\| \leq  M\rme^{-at}\big(|x^0|_p + \|\A\|\big) + \gamma(|u|_{q,\infty}),%
\end{equation*}
and thus%
\begin{equation*}
  |\phi(t,x^0,u)|_p  \leq  M\rme^{-at}|x^0|_p  + \gamma(|u|_{q,\infty}) + (1+ M\rme^{-at}) \|\A\|,%
\end{equation*}
$t \geq 0$.
Consequently,~\eqref{eq_interconnection} is eISS with respect to the origin, but with the ``offset'' $(1+ M\rme^{-at}) \|\A\|$.
This property can be called \emph{exponential input-to-state practical stability}. See~\cite{Mir19a} for more on this property for infinite-dimensional systems.
\end{remark}

For any function $V:X\rightarrow \R$, which is continuous on $X \backslash \A$, we define the \emph{orbital derivative} at $x \in X\backslash \A$ for the external input $u \in \UC$ by
$\rmD^+ V_u(x) := \rmD^+ V(\phi(t,x,u))_{|t=0}$,
where the right-hand side is the right upper Dini derivative of the function $t \mapsto V(\phi(t,x,u))$, evaluated at $t=0$.%

Exponential input-to-state stability is implied by the existence of an exponential ISS Lyapunov function, which we define in a dissipative form as follows.%

\begin{definition}\label{def:eISS-Lyapunov-Function}
Let a nonempty closed set $\A \subset X $ be given. A function $V:X \rightarrow \R_+$, which is continuous on $X \backslash \A$, is called an \emph{eISS Lyapunov function for $\Sigma$ w.r.t.~$\A$} if there exist constants $\underline{\omega},\overline{\omega},b,\kappa>0$ and $\gamma \in \KC_{\infty}$ such that%
\begin{subequations}\label{eq_eiss_lyap_props}
\begin{align}
  \underline{\omega} |x|_\A^b &\leq V(x) \leq \overline{\omega} |x|_\A^b\qquad \forall x\in X,  \label{eq_eiss_lyap_props-1}\\
  \rmD^+ V_u(x) &\leq -\kappa V(x) + \gamma(|u|_{q,\infty}) \quad \forall x\in X\backslash \A, \forall u \in \UC. \label{eq_eiss_lyap_props-2}%
\end{align}
\end{subequations}
The function $\gamma$ is sometimes called a \emph{Lyapunov gain}.
\end{definition}

\begin{proposition}\label{prop_eiss}
If there exists an eISS Lyapunov function for $\Sigma$ w.r.t. $\A$, then $\Sigma$ is eISS w.r.t. $\A$.%
\end{proposition}
The proof follows similar steps as those in the proof of Proposition 4.4 in~\cite{KMS19}.

%

\section{The Gain operator and its properties}\label{sec:The-gain-operator-and-its-properties}

Our main objective is to develop conditions for input-to-state stability of the interconnection of countably many subsystems \eqref{eq_ith_subsystem}, depending on the ISS properties of the subsystems and the interconnection structure. Throughout this section, we assume that the infinite interconnection $\Sigma$ is well-posed with state space $X = \ell^p(\N,(n_i))$ and external input space $U = \ell^q(\N,(m_i))$ for some $p,q \in [1,\infty)$.

\subsection{Assumptions on the subsystems}

We assume that each subsystem $\Sigma_i$, given by \eqref{eq_ith_subsystem}, is exponentially ISS w.r.t.\ a closed set $\A_i$ and there exist continuous eISS Lyapunov functions w.r.t.~$\A_i$ with linear gains for all $\Sigma_i$. The following assumption formulates the eISS property for the subsystems.%

\begin{assumption}\label{ass_vi_existence}
For each $i\in\N$ there is a nonempty closed set $\A_i \subset \R^{n_i}$ and a continuous function $V_i:\R^{n_i} \rightarrow \R_+$, satisfying for certain $p,q\in[1,\infty)$ the following properties.%
\begin{itemize}
\item There are constants $\underline{\alpha}_i,\overline{\alpha}_i>0$ so that for all $x_i \in \R^{n_i}$%
\begin{equation}\label{eq_viest}
  \underline{\alpha}_i|x_i|^p_{\A_i} \leq V_i(x_i) \leq \overline{\alpha}_i|x_i|^p_{\A_i}.%
\end{equation}
\item There are constants $\lambda_i$, $\gamma_{ij}$ ($j \in I_i$), $\gamma_{iu}>0$ so that the following holds: for each $x_i\in\R^{n_i}\backslash \A_i$,  $u_i \in L^{\infty}(\R_+,\R^{m_i})$, each internal input $\bar{x} \in C^0(\R_+,X)$ and for almost all $t$ in the maximal interval of existence of $\phi_i(t) := \phi_i(t,x_i,(\bar{x},u_i))$ one has%
\begin{align}\label{eq_nablaviest}									
\begin{split}
\hspace{-4mm}\rmD^+ (V_i \circ \phi_i)(t) \leq -&\lambda_i V_i(\phi_i(t)) + \sum_{j \in I_i} \gamma_{ij}V_j(x_j(t)) \\
&+ \gamma_{iu}|u_i(t)|^q,%
\end{split}
\end{align}
where we denote the components of $\bar{x}$ by $x_j(\cdot)$.\footnote{At this point, the right-hand side of \eqref{eq_nablaviest} is not necessarily finite. However, this requirement is not needed here.}%
\item For all $t$ in the maximal interval of the existence of $\phi_i$ one has%
\begin{equation*}
  \rmD_+(V_i \circ \phi_i)(t) < \infty.%
\end{equation*}
\end{itemize}
\end{assumption}

We furthermore assume that the following uniformity conditions hold for the constants introduced above.%

\begin{assumption}\label{ass_external_gains}
\begin{enumerate}
\item[(a)] There are constants $\underline{\alpha},\overline{\alpha}>0$ so that for all $i\in \N$
\begin{equation}\label{eq_uniformity_alpha}
  \underline{\alpha} \leq \underline{\alpha}_i \leq \overline{\alpha}_i \leq \overline{\alpha}.%
\end{equation}
\item[(b)] There is a constant $\underline{\lambda}>0$ so that for all $i\in\N$
\begin{equation}\label{eq_uniformity_lambda}
  \underline{\lambda} \leq \lambda_i.%
\end{equation}
\item[(c)] There is a constant $\overline{\gamma}_u>0$ so that for all $i\in\N$%
\begin{equation}\label{eq_uniformity_gammaiu}
  \gamma_{iu} \leq \overline{\gamma}_u.%
\end{equation}
\end{enumerate}
\end{assumption}

In order to formulate a small-gain condition, we further introduce the following infinite nonnegative matrices by collecting the coefficients from \eqref{eq_nablaviest}
\begin{equation*}
  \Lambda := \diag(\lambda_1,\lambda_2,\lambda_3,\ldots), \quad \Gamma := (\gamma_{ij})_{i,j\in\N},%
\end{equation*}
where we put $\gamma_{ij} := 0$ whenever $j \notin I_i$. We also introduce the infinite matrix%
\begin{equation}\label{eq:operator-A}
  \Psi := \Lambda^{-1}\Gamma = (\psi_{ij})_{i,j\in\N},\quad \psi_{ij} = \frac{\gamma_{ij}}{\lambda_i}.
\end{equation}

Under an appropriate boundedness assumption, the matrix $\Psi$ acts as a linear operator on $\ell^1$ by $
  (\Psi x)_i = \sum_{j=1}^{\infty} \psi_{ij} x_j \mbox{\quad for all\ } i \in \N$.

We call $\Psi:\ell^1 \rightarrow \ell^1$ the \emph{gain operator} associated with the decay rates $\lambda_i$ and coefficients $\gamma_{ij}$.%

We make the following assumption, which is equivalent to $\Gamma$ being a bounded operator from $\ell^1$ to $\ell^1$.%

\begin{assumption}\label{ass_A_bounded}
The matrix $\Gamma = (\gamma_{ij})$ satisfies%
\begin{equation}\label{eq_Gamma_bounded}
  \|\Gamma\|_{1,1} = \sup_{j \in \N} \sum_{i=1}^{\infty} \gamma_{ij} < \infty,%
\end{equation}
where the double index on the left-hand side indicates that we consider the operator norm induced by the $\ell^1$-norm both on the domain and codomain of the operator $\Gamma$ (the formula for the norm of $\Gamma$ can be obtained by standard computations as in the case of finite matrices).
\end{assumption}

Clearly, under Assumptions~\ref{ass_A_bounded} and~\ref{ass_external_gains}(b), the gain operator $\Psi$ is bounded (see also \cite[Lem.~V.7]{KMS19}). Moreover, clearly $\Psi$ is a positive operator with respect to the standard positive cone $\ell^1_+:=\{x=(x_1,x_2,\ldots) \in\ell^1:x_i\geq 0,\ \forall i\in\N\}$ in $\ell^1$. Also recall from \cite[Lem.~V.10]{KMS19} the following lemma which uses positive operator theory to deduce the existence of a positive vector $\mu$ that can be used to construct an eISS Lyapunov function for the interconnected system.%

\begin{lemma}\label{lem_smallgain}
Assume that $r(\Psi) < 1$ and that there exists a constant $\overline{\lambda}>0$ such that $\lambda_i \leq \overline{\lambda}$ for all $i \in \N$. Then the following statements hold:%
\begin{itemize}
\item[(i)] There exist a vector $\mu = (\mu_i)_{i\in\N} \in \inner\, \ell^{\infty}_+$ and a constant $\lambda_{\infty}>0$ so that%
\begin{equation}\label{eq_20}
  \frac{[\mu\trn(-\Lambda + \Gamma)]_i}{\mu_i} \leq -\lambda_{\infty} \mbox{\quad for all\ } i \in \N.%
\end{equation}
\item[(ii)] For every $\rho>0$ we can choose the vector $\mu$ and the constant $\lambda_{\infty}$ so that%
\begin{equation}\label{eq_lambdainftyest}
  \lambda_{\infty} \geq (1 - r(\Psi))\underline{\lambda} - \rho.%
\end{equation}
\end{itemize}
\end{lemma}

\section{Small-Gain Theorem}\label{sec:Small-gain-theorem}

In this section, we prove that the interconnected system $\Sigma$ is exponentially ISS under the given assumptions, provided that the spectral radius of the gain operator satisfies $r(\Psi)<1$. By Proposition \ref{prop_eiss}, our objective is reduced to finding an eISS Lyapunov function for the interconnection $\Sigma$, which is accomplished by the following \emph{small-gain theorem}, which is the main result of the paper.

\begin{theorem}\label{MT}
Consider the infinite interconnection $\Sigma$, composed of the subsystems $\Sigma_i$, $i\in\N$, with fixed $p,q \in [1,\infty)$. Suppose that the following hold.%
\begin{enumerate}
\item[(i)] $\Sigma$ is well-posed as a system with state space $X = \ell^p(\N,(n_i))$, space of input values $U = \ell^q(\N,(m_i))$, and the external input space $\UC$, as defined in \eqref{eq:Input-space}.%
\item[(ii)] Each $\Sigma_i$ admits a continuous eISS Lyapunov function $V_i$ w.r.t.~a nonempty closed set $\A_i \subset \R^{n_i}$ so that Assumptions \ref{ass_vi_existence} and \ref{ass_external_gains} are satisfied.%
\item[(iii)] The operator $\Gamma:\ell^1 \rightarrow \ell^1$ is bounded, i.e., Assumption \ref{ass_A_bounded} holds.%
\item[(iv)] The spectral radius of $\Psi$ satisfies $r(\Psi) < 1$.%
\end{enumerate}
Consider the set $\A := X \cap (\A_1 \tm \A_2 \tm \ldots)$. Then $\Sigma$ admits an eISS Lyapunov function w.r.t.~$\A$ of the form%
\begin{equation}\label{eq:Lyapunov-function-construction}
  V(x) = \sum_{i=1}^{\infty} \mu_i V_i(x_i),\quad V:X \rightarrow \R_+%
\end{equation}
for some $\mu = (\mu_i)_{i\in\N}\in \ell^{\infty}$ satisfying $\underline{\mu} \leq \mu_i \leq \overline{\mu}$ with constants $\underline{\mu},\overline{\mu}>0$. In particular, the function $V$ has the following properties.%
\begin{enumerate}
\item[(a)] $V$ is continuous on $X \backslash \A$.%
\item[(b)] There is a $\lambda_\infty>0$ so that for all $x^0 \in X \backslash \A$ and $u \in \UC$%
\begin{equation*}
  \rmD^+ V_u(x^0) \leq -\lambda_{\infty} V(x^0) + \overline{\mu}\,\overline{\gamma}_u|u|_{q,\infty}^q.%
\end{equation*}
\item[(c)] For all $x \in X$ the following inequalities hold:%
\begin{equation}\label{eq:Coercivity-bound-for-V}
  \underline{\mu}\underline{\alpha}|x|_\A^p \leq  V(x) \leq \overline{\mu}\,\overline{\alpha}|x|_\A^p.%
\end{equation}
\end{enumerate}
In particular, $\Sigma$ is eISS w.r.t.~$\A$.%
\end{theorem}

\begin{IEEEproof}
The proof is almost identical to the proof of \cite[Thm.~VI.1]{KMS19}. Hence, we only comment on the differences. An eISS Lyapunov function for $\Sigma$ is defined as in \eqref{eq:Lyapunov-function-construction} with $\mu \in \ell^{\infty}$ given by Lemma \ref{lem_smallgain}. It is well-defined because%
\begin{equation*}
  0 \leq V(x) \leq \sum_{i=1}^{\infty} \mu_i \overline{\alpha}_i |x_i|_{\A_i}^p \leq \overline{\alpha}|\mu|_{\infty} |x|_{\A}^p < \infty.%
\end{equation*}
This shows also the upper bound for \eqref{eq:Coercivity-bound-for-V}. The lower bound for \eqref{eq:Coercivity-bound-for-V} is obtained analogously, and thus inequality \eqref{eq_eiss_lyap_props-1} holds for $V$ (with $b=p$). The proof of continuity of $V$ is almost identical to the case where $\A = \{0\}$, hence we omit the details.%

We now prove the estimate on the orbital derivative under the additional assumption that $\lambda_i \leq \overline{\lambda}$ for all $i\in\N$ with a constant $\overline{\lambda}>0$. Fix an initial state $x^0 \in X\backslash \A$ and an external input $u\in\UC$. We write $\phi(t) = \phi(t,x^0,u)$, $\phi(t) = (\phi_1(t),\phi_2(t),\ldots)$. Then for any $t>0$ (where $\phi(t)$ is defined), we obtain%
\begin{align*}
  \frac{1}{t}\big(V(\phi(t)) - V(x^0)\big) = \frac{1}{t}\sum_{i=1}^{\infty} \mu_i \big[V_i(\phi_i(t)) - V_i(\phi_i(0))\big].%
\end{align*}
Since inequalities \eqref{eq_nablaviest} are valid for almost all positive times, the function in the right-hand side of 
\eqref{eq_nablaviest} is Lebesgue integrable, and since we assume that $\rmD_+(V_i \circ \phi_i)(t)<\infty$ for all $t$, we can proceed using the generalized fundamental theorem of calculus (see \cite[Thm.~9 and p.~42, Rmk.~5.c]{HaT06} or \cite[Thm.~7.3, p.~204]{Sak47}) to%
\begin{align*}
 \frac{1}{t}\big(V(\phi(t)) - V(x^0)\big) &\leq \frac{1}{t}\sum_{i=1}^{\infty}\int_0^t \mu_i \Bigl[ -\lambda_i V_i(\phi_i(s)) \\
	& + \sum_{j\in I_i}\gamma_{ij}V_j(\phi_j(s)) + \gamma_{iu}|u_i(s)|^q\Bigr]\rmd s,%
\end{align*}
where we note that with $\overline{\gamma} := \sup_{i,j}\gamma_{ij}$%
\begin{equation*}
  \sum_{j\in I_i}\gamma_{ij}V_j(\phi_j(s)) \leq \overline{\gamma}\, \overline{\alpha} \sum_{j=1}^{\infty} |\phi_j(s)|_{\A_j}^p < \infty.%
\end{equation*}
We now apply the Fubini-Tonelli theorem in order to interchange the infinite sum and the integral (interpreting the sum as an integral associated with the counting measure on $\N$). To do this, it suffices to prove that the following integral is finite.%
\begin{equation*}
  \int_0^t \sum_{i=1}^{\infty} \Bigl|\mu_i \Bigl[-\lambda_i V_i(\phi_i(s)) + \sum_{j\in I_i}\gamma_{ij}V_j(\phi_j(s)) + \gamma_{iu}|u_i(s)|^q\Bigr]\Bigr| \rmd s.%
\end{equation*}
Using \eqref{eq_viest}, \eqref{eq_uniformity_alpha}, \eqref{eq_uniformity_gammaiu}, and the assumption that $\lambda_i \leq \overline{\lambda}$, we can upper bound the inner term by%
\begin{equation*}
   \overline{\mu}\Bigl[\overline{\lambda}\overline{\alpha}|\phi_i(s)|_{\A_i}^p + \sum_{j \in I_i}\gamma_{ij}\overline{\alpha}|\phi_j(s)|_{\A_j}^p + \overline{\gamma}_u|u_i(s)|^q\Bigr].%
\end{equation*}
By summing the three terms over $i$, one obtains%
\begin{align*}
  \overline{\lambda}\overline{\alpha} \sum_{i=1}^{\infty}|\phi_i(s)|_{\A_i}^p &\leq c_1|\phi(s)|_{\A}^p,\\
	\sum_{i=1}^{\infty}\sum_{j \in I_i}\!\gamma_{ij}\overline{\alpha} |\phi_j(s)|_{\A_j}^p &\!\leq\! c_2  \sum_{j=1}^{\infty} |\phi_j(s)|_{\A_j}^p \sum_{i=1}^{\infty} \gamma_{ij} \!\leq\! c_3  |\phi(s)|_\A^p,\\
  \overline{\gamma}_u\sum_{i=1}^{\infty}|u_i(s)|^q &= c_4  |u(s)|_{q}^q,%
\end{align*}
for some constants $c_1,c_2,c_3,c_4>0$. In the inequality for the middle term, we use the boundedness assumption on the operator $\Gamma$. Hence,%
\begin{align*}
 & \int_0^t \!\sum_{i=1}^{\infty} \Bigl|\mu_i \Bigl[-\lambda_i V_i(\phi_i(s)) + \sum_{j\in I_i}\gamma_{ij}V_j(\phi_j(s)) \!\!+\!\! \gamma_{iu}|u_i(s)|^q\Bigr]\Bigr| \rmd s\\
 &\quad \leq c\int_0^t \left(|\phi(s)|_\A^p + |u(s)|_{q}^q\right) \rmd s < \infty,%
\end{align*}
for some constant $c>0$, where we use the fact that the integrand in the last term is essentially bounded ($s \mapsto |\phi(s)|_\A^p$ is continuous and $s \mapsto |u(s)|_{q}^q$ is essentially bounded).%

Using the notation $
  V_{\mathrm{vec}}(\phi(s)) := (V_1(\phi_1(s)),V_2(\phi_2(s)),\ldots)\trn$
and applying the Fubini-Tonelli theorem, one can conclude that%
\begin{align*}
&\frac{1}{t}\big(V(\phi(t)) - V(x^0)\big) \\
&\leq \frac{1}{t}\!\!\int_0^t \!\!\sum_{i=1}^{\infty} \!\mu_i\! \Bigl[-\lambda_i V_i(\phi_i(s)) \!+\!\! \sum_{j\in I_i}\!\!\gamma_{ij}V_j(\phi_j(s)) \!+\! \gamma_{iu}|u_i(s)|^q\Bigr]\rmd s\\
	&= \frac{1}{t} \int_0^t \Bigl[ \mu\trn(-\Lambda + \Gamma)V_{\mathrm{vec}}(\phi(s)) + \sum_{i=1}^{\infty} \mu_i \gamma_{iu} |u_i(s)|^q\Bigr]\rmd s\\
	&\leq \frac{1}{t}\int_0^t \Bigl[-\lambda_{\infty} V(\phi(s)) + \overline{\mu}\,\overline{\gamma}_u|u|_{q,\infty}^q \Bigr]\rmd s\\
	&= \frac{1}{t}\int_0^t -\lambda_{\infty} V(\phi(s))\, \rmd s + \overline{\mu}\,\overline{\gamma}_u|u|_{q,\infty}^q,%
\end{align*}
where we use \eqref{eq_20} to show the second inequality above. Since $s \mapsto V(\phi(s))$ is continuous, one obtains%
\begin{align*}
  \rmD^+ V_u(x^0) &= \limsup_{t \rightarrow 0+}\frac{1}{t}\left(V(\phi(t)) - V(x^0)\right) \nonumber\\
  &\leq -\lambda_{\infty} V(x^0) + \overline{\mu}\,\overline{\gamma}_u|u|_{q,\infty}^q.%
\end{align*}
Hence, \eqref{eq_eiss_lyap_props-2} holds for $V$ with $\kappa = \lambda_{\infty}$ and $\gamma(r) = \overline{\mu}\,\overline{\gamma}_ur^q$.%

The rest of the proof is identical to the proof of \cite[Thm.~VI.1]{KMS19}, hence we omit the remaining steps.%
\end{IEEEproof}

\section{Applications}\label{sec:Applications}

In this section, we study three applications: stability analysis of time-varying interconnections, dynamic average consensus and the design of distributed observers for infinite networks.%

\subsection{Time-varying interconnected systems}

Although our main result only considers time-invariant systems, it can also be applied to time-varying systems by transforming a time-varying system into a time-invariant one of the form~\eqref{eq_interconnection}. To see this, consider the time-varying system%
\begin{equation}\label{eq:31}
  \dot x = f(t,x,u),%
\end{equation}
where $x \in X$, $u \in U$ and $f \colon \R \times X \times U \to X$ is continuous with $f(t,0,0) = 0$ for all $t \in \R$.
Using the same arguments as those for well-posedness of the network~\eqref{eq_interconnection}, we assume that the state space $X$ and the input space $U$ are chosen as $X = \ell^p(\N,(n_i))$ and $U = \ell^q(\N,(m_i))$, respectively, for fixed $p,q \in [1,\infty)$.
The same class of admissible control functions as in~\eqref{eq:Input-space}	is considered here.

We assume that unique solutions exist for all initial times, initial states and admissible inputs. For any initial time $t^0 \in \R$, initial value $x^0 \in X$ and input $u \in \UC$, the corresponding solution of system~\eqref{eq:31} is denoted by $\phi(\cdot,t^0,x^0,u)$.%


\begin{definition}\label{def:UeISS}
The system~\eqref{eq:31} is called \emph{uniformly exponentially input-to-state stable (UeISS)} if it is forward complete and there are constants $a,M>0$, independent of $t^0$, and $\gamma \in \KC$ such that for any initial time $t^0 \in \R$, initial state $x^0 \in X$ and external input $u\in\UC$ the corresponding solution of~\eqref{eq:31} satisfies for all $t \geq t^0$%
\begin{equation*}
|\phi(t,t^0,x^0,u)|_p \leq M\rme^{-a(t-t^0)}|x^0|_p + \gamma(|u(t^0+\cdot)|_{q,\infty}). %
\end{equation*}
\end{definition}


By adding a ``clock'', one can (see \cite{Teel.2000}) transform~\eqref{eq:31} into%
\begin{align}\label{eq:32}
\begin{split}
  \dot y &= 1, \\
  \dot z &= f(y,z,u),%
\end{split}
\end{align}
where $y \in \R$, $z \in X$, $u \in U$. We equip $\R$ with an arbitrary norm $|\cdot|$ and turn $\R\tm X$ into an $\ell^p$ space by putting%
\begin{equation*}
  |(y,z)|_p := (|y|^p + |z|_p^p)^{1/p}.%
\end{equation*}

Denoting the transition map of \eqref{eq:32} by $\tilde{\phi} = \tilde{\phi}(t,(y,z),u)$, and its $z$-component by $\tilde{\phi}_2$, we see that the following holds:%
\begin{equation}\label{eq:Time-variant-and-time-invariant}
  \phi(t,t^0,x,u) = \tilde{\phi}_2(t-t^0,(t^0,x),u(t^0+\cdot)) \mbox{\quad for all\ } t\geq t^0.%
\end{equation}
The stability properties of \eqref{eq:31} and \eqref{eq:32} are related in the following way:%

\begin{proposition}\label{prop:UeISS criterion} 
The system \eqref{eq:31} is UeISS if and only if \eqref{eq:32} is eISS with respect to the closed set $\A = \{ (y,z) \in \R \times X : z = 0 \} = \R \times \{0\}$.
\end{proposition}

Assume that the system~\eqref{eq:31} can be decomposed into infinitely many interconnected subsystems%
\begin{align}\label{eq:33}
  \dot x_i = f_i(t,x_i,\bar{x},u_i),\quad i\in\N,%
\end{align}
with $t\in\R$, $x_i \in \R^{n_i}$, $\bar{x} \in X$ and $u_i \in \R^{m_i}$. Also, let $f_i \colon \R \tm \R^{n_i} \tm X \tm \R^{m_i} \to \R^{n_i}$ be continuous with $f_i(t,0,0,0) = 0$ for all $t \in \R$.%

With each of the systems \eqref{eq:33} we associate a time-invariant system by%
\begin{align}\label{eq:33-trans}
  \dot z_i = \tilde{f}_i(z_i,(y,\bar{z}),u_i):= f_i(y,z_i,\bar{z},u_i),%
\end{align}
where the time $t$ now becomes an additional internal input $y$. 

Define $\A_0 := \R$ and $\A_i := \{0\} \subset \R^{n_i}$ for all $i \geq 1$. Aggregating all subsystems~\eqref{eq:33-trans}, $i\in\N$, and adding the clock 
$\dot{y} = 1$
as the $0$th subsystem, we obtain an infinite network of the form~\eqref{eq:32}, modeled on the state space $\ell^p(\N_0,(n_i))$ with $n_0 := 1$.%

To enable the stability analysis of the composite system, we make the following assumption.%

\begin{assumption}\label{A:02}
For each $i\in\N$ there exists a continuous function $V_i:\R^{n_i} \rightarrow \R_+$, satisfying for certain $p,q\in[1,\infty)$ the following properties.%
\begin{itemize}
\item There are constants $\underline{\alpha}_i,\overline{\alpha}_i>0$ so that for all $z_i \in \R^{n_i}$%
\begin{equation}\label{eq:viest-time-varying}
  \underline{\alpha}_i|z_i|^p \leq V_i(z_i) \leq \overline{\alpha}_i |z_i|^p.%
\end{equation}
\item There are constants $\lambda_i,\gamma_{ij},\gamma_{iu}>0$ so that the following holds: for each $z_i \in \R^{n_i}$, $u_i \in L^{\infty}(\R_+,\R^{m_i})$ and each internal input $(y,\bar{z}) \in C^0(\R_+,\R \tm X)$ and for almost all $t$ in the maximal interval of existence of $\phi_i(t) := \phi_i(t,z_i,(y,\bar{z},u_i))$ one has%
\begin{align}\label{eq:nablaviest-time-varying}
  \rmD^+ (V_i \circ \phi_i)(t) \leq &-\lambda_i V_i(\phi_i(t)) + \sum_{j \in I_i} \gamma_{ij}V_j(z_j(t)) \nonumber\\
  & + \gamma_{iu}|u_i(t)|^q,%
\end{align}
where we denote the components of $\bar{z}$ by $z_j(\cdot)$.%
\item For all $t$ in the maximal interval of the existence of $\phi_i$ one has $\rmD_+(V_i \circ \phi_i)(t) < \infty$.
\end{itemize}
\end{assumption}

Note that due to the inequalities \eqref{eq_viest} and $\A_0 = \R$, we necessarily have $V_0 = 0$ for the eISS Lyapunov function of the $0$th subsystem (the clock).
Furthermore, we can choose $\lambda_0$ as an arbitrary positive number and $\gamma_{0j}:=0$ for all $j\in\N$.

It follows from Theorem~\ref{MT} that under Assumption~\ref{A:02}, the infinite network of systems~\eqref{eq:33} is uniformly exponentially ISS.
This is summarized by the following corollary of Theorem~\ref{MT}.%

\begin{corollary}\label{C:01}
Consider networks~\eqref{eq:31} and~\eqref{eq:32} and suppose the following:%
\begin{itemize}
\item[(i)] Assumption~\ref{A:02} holds.
\item[(ii)] The constants in Assumption~\ref{A:02} are uniformly bounded as in Assumption \ref{ass_external_gains}.
\item[(iii)] The operator $\Gamma:\ell^1 \rightarrow \ell^1$ is bounded, i.e., Assumption \ref{ass_A_bounded} holds.%
\item[(iv)] The spectral radius of $\Psi$ satisfies $r(\Psi) < 1$.%
\end{itemize}
Then the composite system~\eqref{eq:31} is uniformly eISS.
\end{corollary}

\subsection{Dynamic average consensus}\label{sec:consensus}

Let $G := (\VC,\EC)$ be an undirected graph with the set of nodes $\VC = \N$ and the set of edges $\EC \subseteq \VC \tm \VC$. An edge $(i,j)$ in an undirected infinite graph denotes that nodes $j$ and $i$ exchange information bidirectionally. Node $j$ is an \emph{input neighbor} of node $i$ if $(j,i) \in \EC$. We assume that each agent can only communicate with a finite number of other agents, known as neighbors. Let $\NC_i = \{ j | (j,i) \in \EC \}$ denote the set of the input neighbors of node $i$.%


Let $x_i \in \R^n$ denote the state of node $i\in\VC$. Let each node of $G$ be a (dynamic) agent with dynamics%
\begin{equation}\label{eq_acin_ith_subsystem}
  \Sigma_i:\quad \dot{x}_i = f_i(x_i) + Bu_i,\quad i \in \N,%
\end{equation}
where $u_i \in \R^m$ is the control input, the continuous function $f_i \colon \R^n \to \R^n$ represents the dynamics of each uncoupled node, and $B \in \R^{n \tm m}$.
We model the interconnection $\Sigma$ of these systems on the state space $X := \ell^{\infty}(\N,n)$ with the external input space $U := \ell^{\infty}(\N,m)$ and assume well-posedness for the class of controls $\UC$ as defined before.%

The dynamics in~\eqref{eq_acin_ith_subsystem} do not directly depend on the neighbors' states. But these states might enter the input, i.e., we can define a control law $u_i = q_i(x_i,\ol x_i)$, where $q_i$ is a continuous function on $\R^n\tm\R^{N_i}$, $N_i := |\NC_i|n$, and $\ol x_i \in \R^{N_i}$ is the augmented vector of the states of the neighbors. The aim is to establish control laws, which asymptotically lead to consensus of the agents defined as follows. The agents of the network have reached \emph{consensus} if and only if $x_i = x_j$ for all $i,j \in \VC$. A corresponding state value is called a \emph{consensus point}.%

In several applications of distributed cooperative control, the problem of interest can be formulated as a so-called \emph{dynamic average consensus problem} in which a group of agents cooperates to track a weighted average of locally available time-varying reference signals.
To define a meaningful average of infinitely many quantities, we choose a sequence $(\alpha_i)_{i\in\N}$ of positive real numbers satisfying $\sum_{i=1}^{\infty}\alpha_i = 1$.
One can interpret this sequence as a probability distribution on $\N$.
It is of particular interest to track the following weighted average:
\begin{equation}\label{eq:average}
  x_a := \sum_{i=1}^{\infty} \alpha_i x_i.%
\end{equation}
We observe that for every $x \in X$ we have $
  |x_a| \leq \sum_{i=1}^{\infty} \alpha_i |x_i| \leq |x|_{\infty} < \infty$.


The interconnections of the nodes, which are produced by the control law $q_i$, depend on the strength of the coupling and on the state variables of the nodes. Here we consider the most popular type of coupling which is known as \emph{diffusive coupling}~\cite{Ren.2007}. We assume that the coupling between the $i$th and $j$th agents is defined as a weighted difference, i.e., $a_{ij}(x_i - x_j)$.
Therefore, the control input $u_i$ is given by%
\begin{align}\label{eq:diffusive}
  u_i := -\sigma \sum_{j\in\NC_i}\alpha_j a_{ij}(x_i - x_j),
\end{align}
where $\sigma>0$ denotes the coupling gain between the agents and the interconnections weights $a_{ij}$ satisfy%
\begin{equation}\label{eq:Symmetry-coupling-gains}
  a_{ij} = a_{ji} > 0 ,\quad i,j \in \N \mbox{\quad and \quad} \sup_{i,j} a_{ij} = 1.%
\end{equation}
We assume that $a_{ij} = 0$ for $j \in \N \backslash \NC_i$ and we note that $a_{ij} = 0$ reflects the fact that agent $i$ does not communicate with agent $j$.%

We aim to choose the $a_{ij}$'s and $\sigma$ in~\eqref{eq:diffusive} such that $x_i(t) \to x_j(t) \to x_a(t)$ for all $i,j \in \VC$ as $t \to \infty$. \emph{The difficulty of the dynamic average consensus problem is that each agent is normally connected to only few other agents, and therefore $x_a$ is not available to each agent.}%

For the time derivative of the average we have:
\begin{lemma}\label{lem_average_derivative}
Let the following assumptions be satisfied:%
\begin{enumerate}
\item[(a)] \emph{(Uniform local boundedness)} For every compact set $K \subset \R^n$, there is a $C>0$ so that $|f_i(x)| \leq C$ for all $x \in K$ and $i \in \N$.%
\item[(b)] \emph{(Uniform local Lipschitz continuity)} For every compact set $K \subset \R^n$, there is an $L>0$ so that $|f_i(x) - f_i(y)| \leq L|x - y|$ for all $x,y \in K$ and $i\in\N$.%
\end{enumerate}
Consider a solution $\phi(t) = (\phi_i(t))_{i\in\N}$ of $\Sigma$ corresponding to a continuous control input, defined on some interval $I = [0,T)$. Then $x_a(t) = \sum_{i=1}^{\infty} \alpha_i \phi_i(t)$, $x_a:I \rightarrow \R^n$, is a continuously differentiable function and its derivative satisfies%
\begin{equation}\label{eq:estimate-for-xa}
  \dot{x}_a(t) = \sum_{i=1}^{\infty} \alpha_i \dot{\phi}_i(t) \mbox{\quad for all\ } t \in I.%
\end{equation}
\end{lemma}

\begin{remark}
The condition (a) in the proposition above can be relaxed to uniform boundedness at $0$: There is a $C>0$ so that $|f_i(0)| \leq C$ for all $i\in\N$. It is easy to see that together with the Lipschitz condition (b) this implies uniform boundedness on every compact set.%
\end{remark}

From now on, we assume that the conditions (a) and (b) in Lemma \ref{lem_average_derivative} are satisfied.%

Let us define the error by%
\begin{equation*}
  e_i := \alpha_i (x_i - x_a),\quad i \in \N.%
\end{equation*}
The error vector $e := (e_i)_{i\in\N}$ then satisfies%
\begin{align*}
  |e|_1 &= \sum_{i=1}^{\infty} \left|\alpha_i (x_i - x_a)\right| = \sum_{i=1}^{\infty} \Bigl|\alpha_i x_i - \alpha_i \sum_{j=1}^{\infty}\alpha_j x_j\Bigr| \\
		  	&\leq \sum_{i=1}^{\infty} \alpha_i |x_i| + \sum_{i=1}^{\infty} \alpha_i \sum_{j=1}^{\infty} \alpha_j |x_j|
			  \leq |x|_{\infty} + |x|_{\infty} \\
			  &= 2|x|_{\infty} < \infty.%
\end{align*}
Hence, $e \in \ell^1(\N,n)$. The dynamics of the average is%
\begin{equation*}
  \dot{x}_a = \sum_{i=1}^{\infty} \alpha_i \bigl( f_i(\alpha_i^{-1}e_i + x_a) - \sigma B \sum_{j\in\mathcal{N}_i}\alpha_j a_{ij}(\alpha_i^{-1}e_i - \alpha_j^{-1}e_j)\bigr).%
\end{equation*}
Using the symmetry condition in \eqref{eq:Symmetry-coupling-gains}, one can see that the coupling term vanishes, i.e. $\sum_{i=1}^{\infty} \sum_{j=1}^{\infty} \alpha_i \alpha_j a_{ij}(\alpha_i^{-1}e_i - \alpha_j^{-1}e_j) = 0$.
The convergence of all the sums follows from the estimate $
  \sum_{i=1}^N \sum_{j=1}^M |\alpha_j a_{ij} e_i| \leq |e|_1$ for all $N,M\in\N$.
Hence, we end up with%
\begin{equation}\label{eq_average_dyn}
  \dot{x}_a = \hat{f}_0(\hat{x}):=\sum_{i=1}^{\infty} \alpha_i f_i(\alpha_i^{-1}e_i + x_a),%
\end{equation}
where $\hat{x} = (x_a,e_1,e_2,e_3,\ldots)$.%

The dynamics of the the errors $e_i$, $i\in\N$ is given by%
\begin{align}\label{eq_error_dyn}
  \dot{e}_i &= \hat{f}_i(\hat{x}) = \alpha_i \dot{x}_i  - \alpha_i \dot{x}_a \nonumber\\
	&= \alpha_if_i(\alpha_i^{-1}e_i + x_a) - \alpha_i \sigma B \sum_{j\in\mathcal{N}_i}\alpha_ja_{ij}(\alpha_i^{-1}e_i - \alpha_j^{-1}e_j) \nonumber\\
	& \quad - \alpha_i \sum_{j=1}^{\infty}\alpha_j f_j(\alpha_j^{-1}e_j + x_a).%
\end{align}
Let us write $\hat{\Sigma}_i$ for the $i$th subsystem, where we start the enumeration with $i=0$ so that $x_a$ is the state of the $0$th subsystem and $e_i$ that of the $i$th subsystem for all $i \geq 1$.
The state space of the overall system will be taken to be $\hat{X}:=\ell^1(\N_0,n)$.
Since there is no external input, there is no need to specify an input space $U$.%

The following proposition provides sufficient conditions for well-posedness of $\hat{\Sigma}$.%

\begin{proposition}
\label{prop:Well-posedness of the closed-loop consensus problem}
Assume that all $f_i$ have a common global Lipschitz constant $L>0$ and $f_i(0) = 0$ for all $i\in\N$. Then the interconnection $\hat{\Sigma}$ of the systems $\hat{\Sigma}_i$, $i\in\N$, is well-posed as a system with state space $\hat{X} = \ell^1(\N_0,n)$.%
\end{proposition}

Now we study the stabilization of the average and error system $\hat{\Sigma}$ w.r.t.~the closed set $\A := \R^n \tm \{0\} \tm \{0\} \tm \ldots \subset \ell^1(\N_0,n)$. From our main result, Theorem \ref{MT}, we can immediately conclude the following.%

\begin{theorem} \label{thm:consensus}
Consider the interconnection $\hat{\Sigma}$ of the subsystems $\hat{\Sigma}_i$, $i \in \N_0$, and assume that the following assumptions hold:%
\begin{enumerate}
\item[(i)] The system $\hat{\Sigma}$ with state space $\hat{X} = \ell^1(\N_0,n)$ is well-posed.%
\item[(ii)] Each subsystem $\hat{\Sigma}_i$, $i \geq 1$, admits a continuous eISS Lyapunov function $V_i$ (with respect to the trivial set $\{0\}$, i.e., in the usual sense) so that Assumptions \ref{ass_vi_existence} and \ref{ass_external_gains} are satisfied.%
\item[(iii)] The operator $\Gamma:\ell^1 \rightarrow \ell^1$ is bounded, i.e., Assumption \ref{ass_A_bounded} holds.%
\item[(iv)] The spectral radius of $\Psi$ satisfies $r(\Psi) < 1$.%
\end{enumerate}
Then $\hat{\Sigma}$ is eISS w.r.t.~the set $\A = \R^n \tm \{0\} \tm \{0\} \tm \ldots \subset \ell^1(\N_0,n)$, and thus there are $M>0$ and $a>0$ so that
\begin{eqnarray}
  |e(t)|_1 = \sum_{i=1}^\infty \alpha_i|\phi_i(t)-x_a(t)| \leq M e^{-at} |e(0)|_1.
\label{eq:Weighted-consensus-l-1-norm}
\end{eqnarray}
\end{theorem}

Although Theorem~\ref{thm:consensus} explicitly makes no assumption on the connectedness of the associated graph $G$~\cite{Ren.2007}, the verification of the conditions in the theorem often asks for the connectedness of $G$. We note that in some trivial cases; e.g. if all agents $\Sigma_i$ are linear and individually asymptotically stable, all the conditions will be trivially fulfilled even without making assumptions on connectedness of the agents.%

\begin{remark}
Observe that the local eISS Lyapunov function for the $0$th subsystem has to be chosen identically zero because of \eqref{eq_viest}. Then \eqref{eq_nablaviest} will be trivially satisfied. This is why we do not need any assumption about a local eISS Lyapunov function for the $0$th subsystem. Also observe that the set of neighbors of the $0$th subsystem here is $\N$, while the set of neighbors of the $i$th subsystem is the finite set $\{0\} \cup \NC_i$ for all $i \geq 1$. However, since $V_0 = 0$, the $0$th subsystem does not play any role as a neighbor of the other subsystems.%
\end{remark}

\begin{remark}
Of particular interest in weighted average consensus applications is how to choose the weights $\alpha_i$'s in~\eqref{eq:average}.
As a particular application, in distributed cooperative spectrum sensing, the main objective is to develop distributed protocols for solving the cooperative sensing problem in cognitive radio systems, e.g. see~\cite{Hernandes.2018,Zhang.2015} and references therein.
The weights in this case represent a ratio related to the channel conditions of each agent.
\end{remark}

\subsubsection{Discussion of consensus problem for infinite networks}\label{rem:Discussion on consensus result} 

The estimate \eqref{eq:Weighted-consensus-l-1-norm} shows that the weighted error between the state of the system and the weighted average decreases exponentially. For each $N>0$ define $\alpha^m_{N}:=\min_{i=1}^N\alpha_i$. With this notation we obtain from 
\eqref{eq:Weighted-consensus-l-1-norm} for any $N>0$ that
\begin{eqnarray}
\sum_{i=1}^N |\phi_i(t)-x_a(t)| \leq \frac{M}{\alpha^m_{N}} e^{-at} |e(0)|_1.
\label{eq:Weighted-consensus-l-1-norm-N-modes}
\end{eqnarray}
This means that the trajectories of any finite number of modes exponentially converge to the weighted average consensus point, with the decay rate which does not depend on the number of the agents.

On the other hand, for each mode we have 
\begin{eqnarray}
|\phi_i(t)-x_a(t)| \leq \frac{M}{\alpha_i} e^{-at} |e(0)|_1, \quad i\in\N,
\label{eq:Weighted-consensus-l-1-norm-i-th-mode}
\end{eqnarray}
and since the overshoot $\frac{M}{\alpha_i}$ tends to infinity as $i\to \infty$, 
the estimate \eqref{eq:Weighted-consensus-l-1-norm} does not imply (at least in a straightforward way) even the boundedness of $\sup_{i=1}^\infty |\phi_i(t)-x_a(t)|$.

To see the reasons for the limitations of our approach, note that for the consensus problem the only reasonable state space $X$ for the multi-agent system is $\ell^{\infty}(\N,n)$, as for the systems living in $\ell^{p}(\N,n)$ for $p<\infty$ the uniform convergence to consensus is not possible, unless the consensus point is zero, which reduces the consensus problem to the stabilization problem.

At the same time, our small-gain theorem was derived for the couplings whose state space is in the $\ell_p$ scale, with $p\in(1,+\infty)$.
Therefore, in order to apply the small-gain results to the consensus problems, we scaled the deviations $\phi_i(t)-x_a(t)$ by the coefficients $\alpha_i$ satisfying $\sum_{i=1}^{\infty}\alpha_i = 1$ to ensure that the error $e$ lives in $\ell^{1}(\N,n)$, which enables to use the small-gain theorem developed in this paper.%

\subsection{Distributed observers}\label{sec:distr-observ}

We consider the problem of constructing \emph{distributed observers for networks} of control systems. For simplicity, we set the external inputs $u_i$'s to zero (i.e., $u_i\equiv 0$ for all $i \in\N$) and focus on the network interconnection aspect, rather than discussing the construction of individual local observers.

Our basic assumption is that in a network context, we have local observers of local subsystems. We assume that the states of these \emph{local observers} asymptotically converge to the true state of each subsystem, given perfect knowledge of the true states of neighboring subsystems. Of course such information will be unavailable in practice, and instead each local observer will at best have the state estimates produced by other, neighboring observers available for its operation.%

\subsubsection{The distributed system to be observed}\label{sec:distr-syst-be}


Let the distributed nominal system consist of infinitely many interconnected subsystems%
\begin{equation} \label{eq:27}
  \Sigma_{i}\colon\left\{
  \begin{aligned}
    \dot x_{i} & = f_{i}(x_{i}, \ol x_{i})\\
    y_{i} & = h_{i}(x_{i}, \ol x_{i})\\    
  \end{aligned}
  \right.,\quad i \in \N .
\end{equation}
While $x_i \in \R^{n_i}$ is the state of the system $\Sigma_i$, the quantity $y_i \in \R^{p_i}$ (for some $p_i \in \N$) is the output that can be measured locally and serves as an input for a state observer. 
We denote by $\ol{x}_{i}$ the vector composed of the state variables $x_j$, $j \in I_i$.
Although our general setting allows each subsystem to directly interact with infinitely many other subsystems, in distributed sensing normally each subsystem is only connected to a finite number of subsystems.
Therefore, the set $I_i$ is assumed to be bounded in this application.
To make this observation as clear as possible, in~\eqref{eq:27}, as opposed to the main body of the paper, we use the notation $\ol x_i$ in place of $\ol x$.
Further we assume that $f_{i}\colon\R^{n_{i}}\tm \R^{N_i}  \to\R^{n_i}$ and $h_i\colon\R^{n_{i}}\tm \R^{N_i} \to\R^{p_i}$
are both continuous, where $N_i := \sum_{j \in I_i}n_j$. 

%
\subsubsection{The structure of the distributed observers}
\label{sec:struct-distr-observ}

It is reasonable to assume that a local observer $\mathcal{O}_{i}$ for a system $\Sigma_{i}$ has access to $y_{i}$ and produces an estimate $\hat{x}_{i}$ of $x_{i}$ for all $t \geq 0$. Moreover, we essentially need to know $x_{j}$ for all $j\in I_{i}$ to reproduce the dynamics~\eqref{eq:27}. Access to this kind of information is unrealistic, so instead we assume that it has access to the outputs $y_{j}$ of neighboring subsystems and/or the estimates $\hat x_{j}$ for $j\in I_i$ produced by neighboring observers.
For more details, one may consult the literature on distributed observation and filtering; see e.g.~\cite{Wang.2018} for distributed observers in which the outputs and the state estimates are exchanged among local observers and~\cite{Olfati-Saber.2007} for those in which only state estimates are shared.%

This means that each \emph{local observer} is represented by%
\begin{equation} \label{eq:42}
  \mathcal{O}_{i}\colon \quad \dot{\hat x}_{i} = \hat f_{i} (\hat x_{i}, y_{i}, \ol y_{i},\ol{\hat x}_{i})%
\end{equation}
for some appropriate continuous function $\hat f_{i}$. Here $\ol y_i$ (resp.~$\ol{\hat x}_{i}$) is composed of the outputs $y_j$ (resp.~state variables $\hat x_j$), $j \in I_i$.%

Necessarily, the observers are coupled in the same directional sense as the original distributed subsystems. Based on the small-gain theorem introduced above, this leads us to a framework for the design of distributed observers that guarantees that an interconnection of local observers exponentially tracks the true system state. Thus we consider the composite system given by%
\begin{subequations}
\label{eq:43}
\begin{align}  %
    \dot x_{i} & =  f_{i}(x_{i}, \ol x_{i}), \quad  y_{i} = h_{i}(x_{i}, \ol x_{i}),\\    
     \dot{\hat x}_i & =  \hat f_{i} (\hat x_{i}, y_{i}, \ol y_{i}, \ol{\hat x}_{i}), \quad	i\in \N.
\end{align}
\end{subequations}

\subsubsection{A consistency framework for the design of distributed observers}

Denote by $\phi_i$ and $\hat \phi_i$ the flow maps of the $x_i$-subsystem and $\hat x_i$-subsystem of \eqref{eq:43}, respectively, and define%
\begin{equation*}
  \A_i := \{ (x_i,\hat x_i) \in \Rn[n_i] \tm \Rn[n_i] : x_i = \hat x_i \},\quad i\in\N.
\end{equation*}
Denote also by $\phi$ and $\hat{\phi}$ the flow maps of $x$-subsystem and $\hat{x}$-subsystem of \eqref{eq:43}, respectively.%

\begin{assumption}\label{ass:Vi_existence-observers}
We assume that the sequence of local observers $\mathcal{O}=(\mathcal{O}_{i})_{i\in \N}$ for $\Sigma=(\Sigma_{i})_{i\in \N}$ is given.
Further, there is $p\in[1,\infty)$ so that for each $i\in\N$ there exists a continuous function $V_i \colon \R^{n_i}\to\R_+$, 
as well as constants $\ol\alpha_{i},\ul\alpha_{i}> 0$ and $\lambda_i, \gamma_{ij} > 0$, $j\in I_i$ such that for all $x_i,\hat x_i \in \R^{n_i}$ the following holds:
\begin{equation}
  \label{eq:45} 
\ul\alpha_{i} \abs{x_i - \hat x_i}^p \leq V_i (x_i,\hat x_i) \leq  \ol\alpha_i \abs{x_i-\hat x_i}^p.
\end{equation}
Furthermore, we assume that dissipative estimates
\begin{align}  
\rmD^+ (V_i \circ (\phi_i,\hat \phi_i)) (t) \leq & -\lambda_i V_i (\phi_i(t),\hat \phi_i(t)) \nonumber\\
&+ \sum_{j \in I_i} \gamma_{ij} V_j (x_j(t),\hat x_j(t)) \label{eq:46}
\end{align}
hold for all $i\in\N$ and for all $t$ in the maximal interval of the existence of $\phi_i$ and $\hat\phi_i$ we have
$\rmD_+(V_i \circ (\phi_i,\hat\phi_i))(t) < \infty$.
\end{assumption}

Following our general framework, we choose the state space for the whole system as $X:=\ell^p(\N,(n_i))$ for $p$ as in \eqref{eq:45}.%

\emph{We would like to derive conditions, which ensure that a network of local observers $\mathcal{O}=(\mathcal{O}_{i})_{i\in \N}$ is a \emph{robust distributed observer} for the whole system $\Sigma$, i.e., the error dynamics of the composite system~\eqref{eq:43} is globally exponentially stable.}

Consider $X\times X$ as a Banach space with the norm $\|(x,y)\|_{X\times X}:=\sqrt{|x|_p^2 + |y|_p^2}$, $(x,y)\in X\times X$ and define%
\begin{equation}\label{eq:Set-A-Observer-design}
  \A := \{ (x,\hat x) \in X \tm X : x = \hat x \}= X\cap \A_1\cap \A_2\cap\ldots.%
\end{equation}

We pose the result of this subsection as a corollary, whose proof is, which is a consequence of Theorem~\ref{MT}, is given in Appendix~\ref{proof:MT-Observers}.%

\begin{theorem}\label{thm:MT-Observers}
Consider the infinite interconnection $\Sigma$, given by equations \eqref{eq:27}, and the corresponding composite system \eqref{eq:43}, with fixed $p\in [1,\infty)$. Let the following hold.%
\begin{enumerate}
\item[(i)] \eqref{eq:43} is well-posed as a system on $X\times X$, with $X = \ell^p(\N,(n_i))$ as a state space of $\Sigma$.%
\item[(ii)] Each $\Sigma_i$ admits a continuous eISS Lyapunov function $V_i$ so that Assumptions \ref{ass:Vi_existence-observers}
 and \ref{ass_external_gains} are satisfied.%
\item[(iii)] The operator $\Gamma:\ell^1 \rightarrow \ell^1$ is bounded, i.e., Assumption \ref{ass_A_bounded} holds.%
\item[(iv)] The spectral radius of $\Psi$ satisfies $r(\Psi) < 1$.%
\end{enumerate}
Then the composite system \eqref{eq:43} admits a Lyapunov function w.r.t.~$\A$ as defined in~\eqref{eq:Set-A-Observer-design} of the form%
\begin{equation}\label{eq:Lyapunov-function-construction-observers}
  V(x,\hat{x}) = \sum_{i=1}^{\infty} \mu_i V_i(x_i,\hat{x}_i),\quad V:X \times X\rightarrow \R_+%
\end{equation}
for some $\mu = (\mu_i)_{i\in\N}\in \ell^{\infty}$ satisfying $\underline{\mu} \leq \mu_i \leq \overline{\mu}$ with some constants $\underline{\mu},\overline{\mu}>0$. In particular, the function $V$ has the following properties.%
\begin{enumerate}
\item[(a)] $V$ is continuous on $(X \times X)\backslash \A$.%
\item[(b)] There is a $\lambda_\infty>0$ so that for all $x^0 \in (X\times X) \backslash \A$%
\begin{equation*}
  \rmD^+ V_u(x^0) \leq -\lambda_{\infty} V(x^0).%
\end{equation*}
\item[(c)] For all $x,\hat{x} \in X$ the following inequalities hold%
\begin{equation}\label{eq:Coercivity-bound-for-V-observers}
  \underline{\mu}\underline{\alpha}|(x,\hat{x})|_\A^p \leq  V(x,\hat{x}) \leq \overline{\mu}\,\overline{\alpha}|(x,\hat{x})|_\A^p.%
\end{equation}
\end{enumerate}
Consequently, the error dynamics of \eqref{eq:43} is globally exponentially stable, i.e., there is $\beta\in\KL$ so that the following holds for all $x,\hat x \in X$ and all $t\geq 0$:%
\begin{equation}\label{eq:eISS-Error-dynamics}
  |\phi(t,x)-\hat \phi(t,\hat x)|_p \leq \beta(|x-\hat{x}|_p,t),%
\end{equation}
which in turn means that $\mathcal{O}=(\mathcal{O}_{i})_{i\in \N}$ is a robust distributed observer for $\Sigma$.
\end{theorem}


\section{Conclusions}\label{sec:Conclusions}

We developed a small-gain theorem ensuing exponential ISS with respect to a closed set for infinite networks.
The small-gain condition was given in terms of the spectral radius representing the coupling between participating subsystems, which can be very efficiently checked for a large class of systems. We illustrated the large applicability of our small-gain theorem by applying it to three different natural/engineered problems including stability of time-varying infinite networks at the origin, weighted average consensus, and distributed state estimation.





\appendix
\vspace{-1mm}


\subsection{Proof of Lemma~\ref{lem:Alternative-A-representation}}

First of all, for any $x_i\in\R^{n_i}$ and $z_i \in \A_i$ it holds that%
\begin{equation*}
  |x_i|_{\A_i} = \inf_{y_i \in \A_i} |x_i-y_i| \leq  |x_i-z_i| \leq |x_i|  +  |z_i|.%
\end{equation*}
As $\A \neq\emptyset$, we can choose $z_i\in\A_i$ so that $z=(z_1,z_2,\ldots)\in \A \subset X$. Now for each $N>0$ we have by using the inequality $\gamma(a+b) \leq \gamma(2a) + \gamma(2b)$, which holds for any $\gamma\in\KC$ and all $a,b\geq 0$, that%
\begin{eqnarray}\label{eq:Technical-estimate-distances-1}
  \sum_{i=1}^N |x_i|_{\A_i}^p \leq \sum_{i=1}^N (|x_i| + |z_i|)^p \leq \sum_{i=1}^N (2^p|x_i|^p + 2^p|z_i|^p).%
\end{eqnarray}
As both $x,z\in X$, the limit $N\to\infty$ of the right-hand side exists, and thus%
\begin{equation}
  \sum_{i=1}^\infty |x_i|_{\A_i}^p \leq 2^p(|x|_p^p + |z|_p^p) < \infty.%
\label{eq:Bound-in_A-series}
\end{equation}
Let us show the second claim. Pick any $x \in X$ and any $\tilde{y} \in \A$. Then for every $\varepsilon>0$ there is $N=N(\varepsilon)$ so that%
\begin{eqnarray}\label{eq:Technical-estimate-distances-2}
  \sum_{i=N+1}^\infty|x_i|^p<\frac{\varepsilon}{2^{p+1}},\qquad \sum_{i=N+1}^\infty|\tilde{y}_i|^p<\frac{\varepsilon}{2^{p+1}}.%
\end{eqnarray}
The following holds:%
\begin{align}\label{eq:Technical-estimate-distances-3}	
\begin{split}			
  |x|_{\A} &= \inf_{y\in \A} \Big(\sum_{i=1}^N |x_i-y_i|^p + \sum_{i=N+1}^\infty |x_i-y_i|^p \Big)^{\frac{1}{p}}\\
	  			 &\leq \inf_{y_i \in \A_i,\ i=1,\ldots,N} \Big(\sum_{i=1}^N |x_i-y_i|^p + \sum_{i=N+1}^\infty |x_i - \tilde{y}_i|^p \Big)^{\frac{1}{p}},
\end{split}
\end{align}
where in the last transition we reduce the set of $y$ over which we take an infimum from $\A$ to $\{(y_1,y_2,\ldots,y_N,\tilde{y}_{N+1}, \tilde{y}_{N+2},\ldots): y_i \in \A_i,i=1,\ldots,N\} \subset \A$.%

Estimating the last term in \eqref{eq:Technical-estimate-distances-3} similarly to \eqref{eq:Technical-estimate-distances-1}, and using 
\eqref{eq:Technical-estimate-distances-2}, we obtain%
\begin{align*}
|x|_{\A} 
&\leq \!\!\inf_{y_i \in \A_i,\ i=1,\ldots,N} \!\!\Big(\sum_{i=1}^N \!|x_i-y_i|^p + 2^p\!\!\! \sum_{i=N+1}^\infty\!\! ( |x_i|^p \!+\!|\tilde{y}_i|^p) \Big)^{\frac{1}{p}}\\
				 &\leq \Big(\sum_{i=1}^N \inf_{y_i \in \A_i}|x_i-y_i|^p + \varepsilon \Big)^{\frac{1}{p}} = \Big(\sum_{i=1}^N |x_i|^p_{\A_i} + \varepsilon \Big)^{\frac{1}{p}}.
\end{align*}
By using \eqref{eq:Bound-in_A-series}, we can estimate the last term by
\begin{equation*}
  |x|_{\A} \leq \Big(\sum_{i=1}^\infty |x_i|^p_{\A_i} + \varepsilon \Big)^{\frac{1}{p}}<\infty.%
\end{equation*}
Now, as $\varepsilon>0$ has been chosen arbitrarily, we can take the limit $\varepsilon\to 0$ to obtain%
\begin{equation}\label{eq:above-estim-distance-lemma}
  |x|_{\A} \leq \Big(\sum_{i=1}^\infty |x_i|^p_{\A_i} \Big)^{\frac{1}{p}}.%
\end{equation}
On the other hand, as taking the infimum over all $x \in \A_1\times \A_2\times\ldots$ gives a value not larger than taking the infimum over $\A$, it holds that%
\begin{align*}
  |x|_{\A} &\geq \inf_{y_i\in \A_i,\ i\in\N} \Big(\sum_{i=1}^\infty |x_i-y_i|^p \Big)^{\frac{1}{p}} \\
  &= \Big(\sum_{i=1}^\infty \inf_{y_i\in \A_i}|x_i-y_i|^p \Big)^{\frac{1}{p}} = \Big(\sum_{i=1}^\infty |x_i|^p_{\A_i} \Big)^{\frac{1}{p}},%
\end{align*}
which together with \eqref{eq:above-estim-distance-lemma} completes the proof of the lemma.


\subsection{Proof of Proposition~\ref{prop:UeISS criterion}}

Let \eqref{eq:32} be eISS w.r.t.~$\A$. Pick any $x^0 \in X$, $u\in\UC$, and any $t^0,t \in\R$ so that $t\geq t^0$. By \eqref{eq:Time-variant-and-time-invariant}, for certain $M,a>0$ and some $\gamma\in\Kinf$ it holds that%
\begin{align*}
  |\phi(t,t^0,x^0,u)|_p &= |\tilde{\phi}_2(t-t^0,(t^0,x^0),u(t^0+\cdot))|_p \\
  &= \big|\tilde{\phi}(t-t^0,(t^0,x^0),u(t^0+\cdot))\big|_{\A}\\
	& \leq M\rme^{-a(t-t^0)}|(t^0,x^0)|_{\A} + \gamma(|u(t^0+\cdot)|_{q,\infty}) \\
	&= M\rme^{-a(t-t^0)}|x^0|_p + \gamma(|u(t^0+\cdot)|_{q,\infty}),
\end{align*}
and \eqref{eq:31} is UeISS. For the other direction of the proof, assume that \eqref{eq:31} is UeISS and pick $u \in \UC$, $(t^0,x^0) \in \R \times X$ and $t \geq 0$. Let $\tilde{u} \in \UC$ be defined by%
\begin{equation*}
  \tilde{u}(t) := \left\{\begin{array}{rl}
	                              0 ,& \mbox{if } t \in (-\infty,t^0],\\
																u(t-t^0) ,& \mbox{if } t > t^0.%
												 \end{array}\right.%
\end{equation*}
Then%
\begin{align*}
&  |\tilde{\phi}(t,(t^0,x^0),u)|_{\A} = |\tilde{\phi}_2(t,(t^0,x^0),u)|_p \\
  &= |\tilde{\phi}_2((t+t^0)-t^0,(t^0,x^0),\tilde{u}(t^0+\cdot))|_p \\
&= |\phi(t+t^0,t^0,x^0,\tilde{u})|_p \leq M\rme^{-at}|x^0|_p + \gamma(|\tilde{u}(t^0+\cdot)|_{q,\infty}) \\
																		 &= M\rme^{-at}|(t^0,x^0)|_{\A} + \gamma(|u|_{q,\infty}),%
\end{align*}
and \eqref{eq:32} is eISS w.r.t.~$\A$.%

\vspace{-2mm}

\subsection{Proof of Lemma~\ref{lem_average_derivative}}\label{sec:proof-Lemma-average}

We apply the measure-theoretic version of the Leibniz rule, where we interpret the infinite sum as the integral with respect to the probability measure on $\N$ defined by the sequence $(\alpha_i)_{i\in\N}$. The crucial assumption that we need to apply this rule is that $|\dot{\phi}_i(t)| \leq \theta_i$ for a sequence $(\theta_i)_{i\in\N}$ with $\sum_{i=1}^{\infty}\alpha_i\theta_i < \infty$. We can estimate%
\begin{align*}
  |\dot{\phi}_i(t)| = |f_i(\phi_i(t)) + Bu_i(t)| \leq |f_i(\phi_i(t))| + \|B\||u_i|_{\infty} 
\end{align*}
for all $t \in I$.
Since we assume that $u(\cdot)$ is an element of $L^{\infty}(\R_+,\ell^{\infty}(\N,m))$, we have%
\begin{align*}
  \sum_{i=1}^{\infty} \alpha_i |u_i|_{\infty} &= \sum_{i=1}^{\infty} \alpha_i \esssup_{t\geq0} |u_i(t)| \leq \sum_{i=1}^{\infty} \alpha_i |u|_{\infty,\infty} \\
  &= |u|_{\infty,\infty}  < \infty.%
\end{align*}
It remains to bound $|f_i(\phi_i(t))|$. The initial condition $x^0 = (x^0_i)_{i\in\N}$ satisfies $|x^0_i| \leq R$ for all $i\in\N$ with a constant $R>0$. Pick any $h \in (0,T)$. As $|f_i(x)| \leq C$ on the closed ball of radius $R$ centered at the origin of $\R^n$, and denoting by $L>0$ a uniform Lipschitz constant of the $f_i$ in some ball of radius sufficiently larger than $R$, we obtain for any $t\in[0,h]$%
\begin{align*} 
&  |\phi_i(t) - x^0_i| \leq \int_0^t |f_i(\phi_i(s)) + Bu_i(s)| \rmd s \\
	                    &\leq \int_0^t |f_i(\phi_i(s))| \rmd s + \|B\| \int_0^t |u_i(s)| \rmd s \\
							        &\leq \|B\| h |u|_{\infty,\infty} \!+\! \int_0^t |f_i(\phi_i(s)) - f_i(\phi_i(0))| \rmd s \! +\! h |f_i(\phi_i(0))| \\
											&\leq \underbrace{(h \|B\| |u|_{\infty,\infty} + h C)}_{=: E} + L \int_0^t |\phi_i(s) - x_i^0| \rmd s.%
\end{align*}
Hence, the Gronwall lemma implies $|\phi_i(t)| \leq |x^0_i| + E\rme^{Lh} \leq R + E\rme^{Lh}$ for all $t \in [0,h]$, and consequently $|f_i(\phi_i(t))|$ is bounded by a constant for all $t\in[0,h]$, which implies the desired integrability.%

\vspace{-2mm}

\subsection{Proof of Proposition~\ref{prop:Well-posedness of the closed-loop consensus problem}}

We verify the assumptions made in~\cite[Thm.~3.2]{KMS19}. We denote by $\hat{f}(\hat{x}):=(\hat{f}_i(\hat{x}))_{i\in\N_0}$ the right-hand side of $\hat{\Sigma}$:
\begin{itemize}
\item[(i)] We show that $\hat{f}(\hat{x}) \in \hat{X}$ for all $\hat{x} \in \hat{X}$. To this end, it suffices to prove that%
\begin{align*}
&\!\! \sum_{i=1}^{\infty}\!|\alpha_if_i(\alpha_i^{-1}e_i + x_a) - \alpha_i \sigma B \!\!\sum_{j\in\mathcal{N}_i}\alpha_ja_{ij}(\alpha_i^{-1}e_i \!-\! \alpha_j^{-1}e_j) \\
  &- \alpha_i \sum_{j=1}^{\infty}\alpha_j f_j(\alpha_j^{-1}e_j + x_a)| < \infty%
\end{align*}
whenever $(x_a,e) \in \hat{X}$. This is a consequence of the following estimates:%
\begin{align*}
&  \sum_{i=1}^{\infty} \alpha_i |f_i(\alpha_i^{-1}e_i + x_a)| \leq \sum_{i=1}^{\infty} \alpha_i L |\alpha_i^{-1}e_i + x_a| \\
&\leq L(|e|_1 + |x_a|),\\
&\sum_{i=1}^{\infty} \Bigl|\alpha_i \sigma B \sum_{j=1}^{\infty} \alpha_j a_{ij}(\alpha_i^{-1}e_i - \alpha_j^{-1}e_j)\Bigr| \\
&\leq \sigma \|B\| \sum_{i=1}^{\infty}\sum_{j=1}^{\infty} \alpha_i\alpha_j  |\alpha_i^{-1}e_i - \alpha_j^{-1}e_j| \\
& \leq \mbox{const} \sum_{i,j=1}^{\infty}\alpha_j |e_i| = \mbox{const} |e|_1,\\
&	\sum_{i=1}^{\infty}\Bigl|\alpha_i \sum_{j=1}^{\infty}\alpha_j f_j(\alpha_j^{-1}e_j + x_a)\Bigr| \leq L(|e|_1 + |x_a|).%
\end{align*}
\item[(ii)] The Lipschitz continuity of the right-hand side is shown as follows:%
\begin{align*}
  & \!\!\!\!\!\!\!\!\!\! \sum_{i=1}^{\infty} \Bigl|\alpha_if_i(\alpha_i^{-1}e_i + x_a) - \alpha_i f_i(\alpha_i^{-1}\tilde{e}_i + \tilde{x}_a)\Bigr| \\
  &\!\!\!\!\!\!\!\!\!\!\leq \sum_{i=1}^{\infty} L (|e_i - \tilde{e}_i| + \alpha_i |x_a - \tilde{x}_a|)\\
	&\!\!\!\!\!\!\!\!\!\!\leq L (|e - \tilde{e}|_1 + |x_a - \tilde{x}_a|) = L |x - \hat{x}|,\\
&	\!\!\!\!\!\!\!\!\!\! \sum_{i=1}^{\infty}\Bigl|\alpha_i \sigma B \sum_{j\in\NC_i}\alpha_ja_{ij}(\alpha_i^{-1}e_i - \alpha_j^{-1}e_j) \\
& \!\!\!\!\!\!\!\!\!\!\qquad - \alpha_i \sigma B \sum_{j\in\NC_i}\alpha_ja_{ij}(\alpha_i^{-1}\tilde{e}_i - \alpha_j^{-1}\tilde{e}_j)\Bigr|\\
& \!\!\!\!\!\!\!\!\!\!\leq \sigma \|B\|\sum_{i=1}^{\infty} \alpha_i\sum_{j=1}^{\infty}\alpha_j|(\alpha_i^{-1}e_i \!-\! \alpha_j^{-1}e_j) \!- \!(\alpha_i^{-1}\tilde{e}_i \!-\! \alpha_j^{-1}\tilde{e}_j)|\\
&\!\!\!\!\!\!\!\!\!\!\leq \mbox{const} \sum_{i=1}^{\infty} \alpha_i \sum_{j=1}^{\infty} \alpha_j (\alpha_i^{-1}|e_i - \tilde{e}_i| + \alpha_j^{-1}|e_j - \tilde{e}_j|)\\
&\!\!\!\!\!\!\!\!\!\!\leq \mbox{const} |e - \tilde{e}|_1.%
\end{align*}
\end{itemize}

\subsection{Proof of Theorem~\ref{thm:MT-Observers}}\label{proof:MT-Observers}

Applying Theorem~\ref{MT}, we obtain that $V$ is an exponential Lyapunov function for the composite system \eqref{eq:43} with respect to the set $\A$.%

The distance of $(x,y)\in X\times X$ to the set $\A$ can be computed as
\begin{align*}
  |(x,y)|_\A &:= \inf_{z \in X}\|(x,y) - (z,z)\|_{X\times X} \\
  & =  \inf_{z \in X} \sqrt{|x-z|_p^2 + |y-z|_p^2}= \frac{1}{\sqrt{2}}|x-y|_p,%
\end{align*}
where the infimum is achieved at $z = \frac{1}{2}(x+y)$. This allows us to represent the norm of the error $e(t,x,\hat{x}):=\phi(t,x)-\hat \phi(t,\hat x)$ of the observer system \eqref{eq:43} as%
\begin{equation}\label{eq:Error-dynamics}
 |e(t,x,\hat{x})|_p :=|\phi(t,x)-\hat \phi(t,\hat x)|_p = \sqrt{2}\big|\big(\phi(t,x),\hat \phi(t,\hat x)\big)\big|_\A.
\end{equation}
Hence, global exponential stability of \eqref{eq:43} w.r.t.~$\A$ implies global exponential stability of the error dynamics (w.r.t.~the $X$-norm).


\end{document}